# Invariance principles for random walks conditioned to stay positive

Francesco Caravenna[a] and Loïc Chaumont[b]

[a]*Dipartimento di Matematica Pura e Applicata, Università di Padova, via Trieste 63, 35121 Padova, Italy.
E-mail: francesco.caravenna@math.unipd.it*
[b]*Laboratoire de Probabilités et Modèles Aléatoires, Université Pierre et Marie Curie, 4 Place Jussieu, 75252 Paris Cedex 05, France. E-mail: chaumont@ccr.jussieu.fr*

**Abstract.** Let $\{S_n\}$ be a random walk in the domain of attraction of a stable law $\mathcal{Y}$, i.e. there exists a sequence of positive real numbers $(a_n)$ such that $S_n/a_n$ converges in law to $\mathcal{Y}$. Our main result is that the rescaled process $(S_{\lfloor nt \rfloor}/a_n, t \geq 0)$, when conditioned to stay positive, converges in law (in the functional sense) towards the corresponding stable Lévy process conditioned to stay positive. Under some additional assumptions, we also prove a related invariance principle for the random walk killed at its first entrance in the negative half-line and conditioned to die at zero.

**Résumé.** Soit $\{S_n\}$ une marche aléatoire dont la loi est dans le domaine d'attraction d'une loi stable $\mathcal{Y}$, i.e. il existe une suite de réels positifs $(a_n)$ telle que $S_n/a_n$ converge en loi vers $\mathcal{Y}$. Nous montrons que le processus renormalisé $(S_{\lfloor nt \rfloor}/a_n, t \geq 0)$, une fois conditionné à rester positif, converge en loi (au sens fonctionnel) vers le processus de Lévy stable de loi $\mathcal{Y}$ conditionné à rester positif. Sous certaines hypothèses supplémentaires, nous montrons un principe d'invariance pour cette marche aléatoire tuée lorsqu'elle quitte la demi-droite positive et conditionnée à mourir en 0.

## 1. Introduction and main results

A well known invariance principle asserts that if the random walk $(S = \{S_n\}, \mathbb{P})$ is in the domain of attraction of a stable law, with norming sequence $(a_n)$, then under $\mathbb{P}$ the rescaled process $\{S_{\lfloor nt \rfloor}/a_n\}_{t \geq 0}$ converges in law as $n \to \infty$ towards the corresponding stable Lévy process, see [27]. Now denote by $(S, \mathbb{P}_y^+)$ the random walk starting from $y \geq 0$ and conditioned to stay always positive (one can make sense of this by means of an $h$-transform, see below or [3]). Then a natural question is whether the rescaled process obtained from $(S, \mathbb{P}_y^+)$ converges in law to the corresponding stable Lévy process conditioned to stay positive in the same sense, as defined in [11].

The main purpose of this paper is to show that the answer to the above question is positive in full generality, i.e. with no extra assumption other than the positivity of the index of the limiting stable law (in order for the conditioning to stay positive to make sense). Under some additional assumptions, we also prove a related invariance principle for the random walk killed at the first time it enters the negative half-line and conditioned to die at zero. Before stating precisely our results, we recall the essentials of the conditioning to stay positive for random walks and Lévy processes.

*1.1. Random walk conditioned to stay positive*

We denote by $\Omega_{\mathrm{RW}} := \mathbb{R}^{\mathbb{Z}^+}$, where $\mathbb{Z}^+ := \{0, 1, \ldots\}$, the space of discrete trajectories and by $S := \{S_n\}_{n \in \mathbb{Z}^+}$ the coordinate process which is defined on this space:

$$\Omega_{\mathrm{RW}} \ni \xi \longmapsto S_n(\xi) := \xi_n.$$

Probability laws on $\Omega_{\mathrm{RW}}$ will be denoted by blackboard symbols.

Let $\mathbb{P}_x$ be the law on $\Omega_{\mathrm{RW}}$ of a random walk started at $x$, that is $\mathbb{P}_x(S_0 = x) = 1$ and under $\mathbb{P}_x$ the variables $\{S_n - S_{n-1}\}_{n \in \mathbb{N} := \{1, 2, \ldots\}}$ are independent and identically distributed. For simplicity, we put $\mathbb{P} := \mathbb{P}_0$. Our basic assumption is that *the random walk oscillates*, i.e. $\limsup_k S_k = +\infty$ and $\liminf_k S_k = -\infty$, $\mathbb{P}$-a.s.

Next we introduce the *strict descending* ladder process, denoted by $(\overline{T}, \overline{H}) = \{(\overline{T}_k, \overline{H}_k)\}_{k \in \mathbb{Z}^+}$, by setting $\overline{T}_0 := 0, \overline{H}_0 := 0$ and

$$\overline{T}_{k+1} := \min\{j > \overline{T}_k \colon -S_j > \overline{H}_k\}, \quad \overline{H}_k := -S_{\overline{T}_k}.$$

Note that under our hypothesis $\overline{T}_k < \infty$, $\mathbb{P}$-a.s. for all $k \in \mathbb{Z}^+$, and that $(\overline{T}, \overline{H})$ is under $\mathbb{P}$ a *bivariate renewal process*, that is a random walk on $\mathbb{R}^2$ with step law supported in the first quadrant. We denote by $V$ the *renewal function* associated to $\overline{H}$, that is the positive nondecreasing right-continuous function defined by

$$V(y) := \sum_{k \geq 0} \mathbb{P}(\overline{H}_k \leq y), \quad y \geq 0. \tag{1.1}$$

Notice that $V(y)$ is the expected number of ladder points in the stripe $[0, \infty) \times [0, y]$. It follows in particular that the function $V(\cdot)$ is *subadditive*.

The only hypothesis that $\limsup_k S_k = +\infty$, $\mathbb{P}$-a.s., entails that the function $V(\cdot)$ is invariant for the semigroup of the random walk killed when it first enters the negative half-line (see Appendix B). Then for $y \geq 0$ we denote by $\mathbb{P}_y^+$ the $h$-transform of this process by $V(\cdot)$. More explicitly, $(S, \mathbb{P}_y^+)$ is the Markov chain whose law is defined for any $N \in \mathbb{N}$ and for any $B \in \sigma(S_1, \ldots, S_N)$ by

$$\mathbb{P}_y^+(B) := \frac{1}{V(y)} \mathbb{E}_y(V(S_N) \mathbf{1}_B \mathbf{1}_{\mathcal{C}_N}), \tag{1.2}$$

where $\mathcal{C}_N := \{S_1 \geq 0, \ldots, S_N \geq 0\}$. We call $\mathbb{P}_y^+$ the law of the random walk starting from $y \geq 0$ and *conditioned to stay positive*. This terminology is justified by the following result which is proved in [3], Theorem 1:

$$\mathbb{P}_y^+ := \lim_{N \to \infty} \mathbb{P}_y(\cdot | \mathcal{C}_N), \quad y \geq 0. \tag{1.3}$$

Note that since in some cases $\mathbb{P}_y^+(\min_{k \geq 0} S_k = 0) > 0$, we should rather call $\mathbb{P}_y^+$ the law of the random walk starting from $y \geq 0$ and conditioned to stay *non-negative*. The reason for which we misuse this term is to fit in with the usual terminology for stable Lévy processes.

We point out that one could also condition the walk to stay strictly positive: this amounts to replacing $\mathcal{C}_N$ by $\mathcal{C}_N^\sim := \{S_1 > 0, \ldots, S_N > 0\}$ and $V(\cdot)$ by $V^\sim(x) := V(x-)$, see Appendix B. The extension of our results to this case is straightforward.

### 1.2. Lévy process conditioned to stay positive

We introduce the space of real-valued càdlàg paths $\Omega := D([0,\infty), \mathbb{R})$ on the real half line $[0,\infty)$, and the corresponding coordinate process $X := \{X_t\}_{t \geq 0}$ defined by

$$\Omega \ni \omega \longmapsto X_t(\omega) := \omega_t.$$

We endow $\Omega$ with the Skorohod topology, and the natural filtration of the process $\{X_t\}_{t \geq 0}$ will be denoted by $\{\mathcal{F}_t\}_{t \geq 0}$. Probability laws on $\Omega$ will be denoted by boldface symbols.

Let $\mathbf{P}_x$ be the law on $\Omega$ of a stable process started at $x$. As in discrete time, we set $\mathbf{P} := \mathbf{P}_0$. Let $\alpha \in (0, 2]$ be the index of $(X, \mathbf{P})$ and $\rho$ be its positivity parameter, i.e. $\mathbf{P}(X_1 \geq 0) = \rho$. When we want to indicate explicitly the parameters, we will write $\mathbf{P}_x^{[\alpha, \rho]}$ instead of $\mathbf{P}_x$. We assume that $\rho \in (0, 1)$ (that is we are excluding subordinators and cosubordinators) and we set $\overline{\rho} := 1 - \rho$. We recall that for $\alpha > 1$ one has the constraint $\rho \in [1 - 1/\alpha, 1/\alpha]$, hence $\alpha \rho \leq 1$ and $\alpha \overline{\rho} \leq 1$ in any case.

We introduce the law on $\Omega$ of the Lévy process starting from $x > 0$ and conditioned to stay positive on $(0, \infty)$, denoted by $\mathbf{P}_x^+$, see [11]. As in discrete time, $\mathbf{P}_x^+$ is an $h$-transform of the Lévy process killed when it first enters the negative half-line, associated to the positive invariant function given by

$$\widetilde{U}(x) := x^{\alpha \overline{\rho}}. \tag{1.4}$$

More precisely, for all $t \geq 0$, $A \in \mathcal{F}_t$ and $x > 0$ we have

$$\mathbf{P}_x^+(A) := \frac{1}{\widetilde{U}(x)} \mathbf{E}_x(\widetilde{U}(X_t) \mathbf{1}_A \mathbf{1}_{\{\underline{X}_t \geq 0\}}), \tag{1.5}$$

where $\underline{X}_t = \inf_{0 \leq s \leq t} X_s$. In analogy with the random walk case, $\widetilde{U}(\cdot)$ is the renewal function of the ladder heights process associated to $-X$, see [2].

We stress that a continuous time counterpart of the convergence (1.3) is valid, see [11], but the analogies between discrete and continuous time break down for $x = 0$. In fact definition (1.5) does not make sense in this case, and indeed 0 is a boundary point of the state space $(0, \infty)$ for the Markov process $(X, \{\mathbf{P}_x^+\}_{x > 0})$. Nevertheless, it has been shown in [12] that it is still possible to construct the law $\mathbf{P}^+ := \mathbf{P}_0^+$ of a càdlàg Markov process with the same semigroup as $(X, \{\mathbf{P}_x^+\}_{x > 0})$ and such that $\mathbf{P}^+(X_0 = 0) = 1$, and we have

$$\mathbf{P}_x^+ \implies \mathbf{P}^+, \quad \text{as } x \downarrow 0,$$

where here and in the sequel $\Rightarrow$ means convergence in law (in particular, when the space is $\Omega$, this convergence is to be understood in the functional sense).

### 1.3. A first invariance principle

The basic assumption underlying this work is that $\mathbb{P}$ is the law on $\Omega_{\mathrm{RW}}$ of a random walk which is attracted to $\mathbf{P}^{[\alpha, \rho]}$, law on $\Omega$ of a stable Lévy process with index $\alpha$ and positivity parameter $\rho \in (0, 1)$. More explicitly, we assume that there exists a positive sequence $(a_n)$ such that as $n \to \infty$

$$S_n / a_n \text{ under } \mathbb{P} \implies X_1 \text{ under } \mathbf{P}^{[\alpha, \rho]}. \tag{1.6}$$

Observe that the hypothesis $\rho \in (0, 1)$ entails that the random walk $(S, \mathbb{P})$ oscillates, so that all the content of Section 1.1 is applicable. In particular, the law $\mathbb{P}_y^+$ is well defined.

Next we define the rescaling map $\phi_N : \Omega_{\mathrm{RW}} \to \Omega$ defined by

$$\Omega_{\mathrm{RW}} \ni \xi \longmapsto (\phi_N(\xi))(t) := \frac{\xi_{\lfloor Nt \rfloor}}{a_N}, \quad t \in [0, \infty). \tag{1.7}$$

For $x \geq 0$ and $y \geq 0$ such that $x = y/a_N$, we define the probability laws

$$\mathbf{P}_x^N := \mathbb{P}_y \circ (\phi_N)^{-1}, \qquad \mathbf{P}_x^{+,N} := \mathbb{P}_y^+ \circ (\phi_N)^{-1}, \tag{1.8}$$

which correspond respectively to the push-forwards of $\mathbb{P}_y$ and $\mathbb{P}_y^+$ by $\phi_N$. As usual, we set $\mathbf{P}^N := \mathbf{P}_0^N$ and $\mathbf{P}^{+,N} := \mathbf{P}_0^{+,N}$. We can now state our main result.

**Theorem 1.1.** *Assume that Eq.* (1.6) *holds true, with $\rho \in (0,1)$, and let $(x_N)$ be a sequence of nonnegative real numbers such that $x_N \to x \geq 0$ as $N \to \infty$. Then one has the following weak convergence on $\Omega$:*

$$\mathbf{P}_{x_N}^{+,N} \implies \mathbf{P}_x^+, \quad N \to \infty. \tag{1.9}$$

The proof of this theorem is first given in the special case when $x_N = x = 0$ for every $N$, see Section 3. The basic idea is to use the absolute continuity between $\mathbf{P}^{+,N}$ (resp. $\mathbf{P}^+$) and the meander of $(X, \mathbf{P}^N)$ (resp. $(X, \mathbf{P})$) and then to apply the weak convergence of these meanders which has been proved by Iglehart [24], Bolthausen [6] and Doney [14].

The proof in the general case is then given in Section 4. The key ingredient is a path decomposition of the Markov chain $(S, \mathbb{P}_y^+)$ at its overall minimum, which is interesting in itself and is presented in Appendix A.

We stress that Theorem 1.1 is valid also in the case $\rho = 1$. The proof of this fact is even simpler than for the case $\rho \in (0,1)$, but it has to be handled separately and we omit it for brevity.

### 1.4. Conditioning to die at zero

Next we consider another Markov chain connected to the positivity constraint: the random walk started at $y \geq 0$, killed at its first entrance in the nonpositive half-line and *conditioned to die at zero*, denoted by $(S, \mathbb{P}_y^\searrow)$. We focus for simplicity on the lattice case, that is we assume that $S_1$ is $\mathbb{Z}$–valued and aperiodic, and we introduce the stopping times $\boldsymbol{\zeta} := \inf\{n \in \mathbb{Z}^+ : S_n = 0\}$ and $T_{(-\infty,0]} := \inf\{n \in \mathbb{Z}^+ : S_n \leq 0\}$. Then the process $(S, \mathbb{P}_y^\searrow)$ is defined as follows: for $y > 0$ we set for $N \in \mathbb{N}$ and $B \in \sigma(S_1, \ldots, S_N)$

$$\mathbb{P}_y^\searrow(B, \boldsymbol{\zeta} > N) := \mathbb{P}_y(B, T_{(-\infty,0]} > N | S_{T_{(-\infty,0]}} \in (-1,0]),$$

and $\mathbb{P}_y^\searrow(S_{\boldsymbol{\zeta}+n} := 0 \; \forall n \geq 0) = 1$, while for $y = 0$ we just set $\mathbb{P}_0^\searrow(S \equiv 0) = 1$.

In an analogous way one could define the Lévy process conditioned to die at zero $(X, \mathbf{P}_x^\searrow)$ (we refer to Section 5.1 for more details) and our purpose is to obtain the corresponding invariance principle. To this aim we introduce the rescaled law of $\mathbb{P}_y^\searrow$ on $\Omega$, by setting for all $x, y \geq 0$ such that $x = y/a_N$

$$\mathbf{P}_x^{\searrow,N} := \mathbb{P}_y^\searrow \circ (\phi_N)^{-1}. \tag{1.10}$$

As it will be shown in Section 5, the process $(S, \mathbb{P}_y^\searrow)$ is an $h$-transform of the random walk $(S, \mathbb{P}_y)$ killed at the first time it enters the nonpositive half-line, corresponding to the excessive function $W(y) := V(\lceil y \rceil) - V(\lceil y \rceil - 1)$, where $V(\cdot)$ is the renewal function defined in (1.1). We will show in Lemma 2.1 that

$$V(x) \sim \frac{x^{\alpha \overline{\rho}}}{L(x)}, \quad x \to \infty, \tag{1.11}$$

for some function $L(\cdot)$ that is slowly varying at infinity. The basic extra-assumption we need to make in order to prove the invariance principle for $(S, \mathbb{P}_y^\searrow)$ is that $W(\cdot)$ satisfies the local form of the above asymptotic relation, namely

$$W(x) \sim \frac{\alpha \overline{\rho}}{L(x)} x^{\alpha \overline{\rho} - 1}, \quad x \to \infty. \tag{1.12}$$

This relation can be viewed as a local renewal theorem for the renewal process $\{\overline{H}_k\}$. Actually, a result of Garsia and Lamperti [21] shows that when $\alpha \overline{\rho} > 1/2$ Eq. (1.12) follows from (1.11), so that we are making no extra-assumption. However in general, i.e. for a generic renewal function $V(\cdot)$, when $\alpha\overline{\rho} \leq 1/2$ Eq. (1.12) is stronger than (1.11). As a matter of fact our setting is very peculiar and it is likely that Eq. (1.12) holds true for all values of $\alpha$ and $\overline{\rho}$, but this remains to be proved.

**Remark 1.2.** *A related open problem concerns the asymptotic behaviour of the probability tail $\mathbb{P}(\overline{H}_1 \geq x)$. In fact by standard Tauberian Theorems [4], Section 8.6.2, we have that for $\alpha\overline{\rho} < 1$ Eq. (1.11) is equivalent to the relation*

$$\mathbb{P}(\overline{H}_1 \geq x) \sim \frac{1}{\Gamma(1+\alpha\overline{\rho})\Gamma(1-\alpha\overline{\rho})} \frac{L(x)}{x^{\alpha\overline{\rho}}}, \quad x \to \infty,$$

*where $\Gamma(\cdot)$ is Euler's Gamma function. The open question is whether the local version of this relation holds true, that is (in the lattice case) whether for $x \in \mathbb{N}$ one has $\mathbb{P}(\overline{H}_1 = x) \sim (\alpha\overline{\rho})x^{-1}\mathbb{P}(\overline{H}_1 \geq x)$ as $x \to \infty$. If this were the case, then Eq. (1.12) would hold true as a consequence of Theorem B in [17].*

We are finally ready to state the invariance principle for the process $(S, \mathbb{P}_y^{\searrow})$. The proof is given in Section 5.

**Theorem 1.3.** *Assume that the law $\mathbb{P}(S_1 \in \mathrm{d}x)$ is supported in $\mathbb{Z}$ and is aperiodic. Assume moreover that Eq. (1.6) holds, with $\rho \in (0,1)$, and that Eq. (1.12) holds true (which happens for instance when $\alpha\overline{\rho} > 1/2$). Let $(x_N)$ be a sequence of nonnegative real numbers that converges towards $x \geq 0$. Then we have the following weak convergence on $\Omega$:*

$$\mathbf{P}_{x_N}^{\searrow, N} \implies \mathbf{P}_x^{\searrow}, \quad N \to \infty.$$

### 1.5. Some motivations and a look at the literature

The study of invariance principles for random walks is a very classical theme in probability theory, cf. [27] and [5]. The extension of these invariance principles to *conditioned* random walks is typically not straightforward: unless a clever representation of the conditioned process can be given, like for instance in [6] and [14] for the meander, quite some technical effort is required, cf. [24] for the meander and [26] for the bridge.

The case of random walks conditioned to stay always positive is analogous: in fact in the specific instance $\alpha = 2, \rho = 1/2$, i.e. when $(S, \mathbb{P})$ is attracted to the Gaussian law, Theorem 1.1 has been proven recently by Bryn-Jones and Doney [7], by combining tightness arguments with a suitable local limit theorem for the conditioned random walk (that has been obtained independently also in [8]). The purpose of our paper is to show that Theorem 1.1 can be proven in a more straightforward way and with lighter techniques, exploiting the absolute continuity with the meander process, see Section 3. Moreover, our proof has the advantage of being completely general, holding for the case of a generic stable law.

Besides the interest they have in their own, invariance principles are important in view of their applications: let us mention two of them that are relevant for our setting. The first application concerns the field of *Lévy trees*. In the stable case, these metric spaces may be obtained as limits of rescaled Galton–Watson trees whose offspring distribution is in the domain of attraction of a stable law. Rather than considering random metric spaces it is often simpler to work with their coding processes. In discrete time, it is known that downward skipfree random walks code the genealogy of Galton–Watson trees and that Lévy processes with no negative jumps code the genealogy of Lévy trees, see for instance Duquesne and Le Gall [18].

A Galton–Watson tree with immigration is obtained by adding independently from a Galton–Watson tree, at each generation $n$, a number $Y_n$ of particles, where the $Y_i$'s are i.i.d. These objects have continuous counterparts, the Lévy trees with immigration, which have been defined by Lambert [25]. Similarly to the standard case, the genealogy of Galton–Watson and Lévy trees with immigration can be coded respectively by random walks and Lévy processes conditioned to stay positive, see Theorem 7 in [25]. In this context, Theorem 1.1 provides a way to prove that an invariance principle still holds in the case of branching processes with immigration.

The second application we present is in the context of random walk models for *polymers and interfaces* (see [22] and [28] for an overview). Consider for instance the following model: for $N \in \mathbb{N}$, $y \in \mathbb{Z}^+$ and $\varepsilon \in \mathbb{R}$ we set

$$\frac{\mathrm{d}\mathbb{P}_{N,y,\varepsilon}}{\mathrm{d}\mathbb{P}_y}(S) := \frac{1}{Z_{N,y,\varepsilon}} \exp\left(\varepsilon \sum_{n=1}^{N} \mathbf{1}_{(S_n=0)}\right) \mathbf{1}_{(S_1 \geq 0, \ldots, S_N \geq 0)}, \tag{1.13}$$

where $(S, \mathbb{P}_y)$ is an aperiodic $\mathbb{Z}$-valued random walk in the domain of attraction of a stable law and $Z_{N,y,\varepsilon}$ is the normalizing constant. The law $\mathbb{P}_{N,y,\varepsilon}$ may be viewed as an effective model for a $(1+1)$-dimensional interface above a hard wall by which the interface is attracted or repelled (depending on the sign of $\varepsilon$). The goal is to find the asymptotic behaviour of the typical paths of $\mathbb{P}_{N,y,\varepsilon}$ in the limit $N \to \infty$ and to study its dependence on $y$ and $\varepsilon$.

The crucial observation is that the measure $\mathbb{P}_{N,y,\varepsilon}$ exhibits a remarkable decoupling between the zero level set $\mathcal{I}_N := \{n \in \{1, \ldots, N\}: S_n = 0\}$ and the excursions of $S$ between two consecutive zeros (we refer to [10] and [9] for more details). In fact, conditionally on $\mathcal{I}_N = (t_1, \ldots, t_k)$, the *bulk excursions* $e_i = \{e_i(n)\}_n := \{S_{t_i+n}\}_{0 \leq n \leq t_{i+1}-t_i}$, for $i = 1, \ldots, k-1$, are independent under $\mathbb{P}_{N,y,\varepsilon}$ and they are distributed like excursions of the free random walk $(S, \mathbb{P}_0)$ conditioned to have a fixed length $(t_{i+1} - t_i)$. Analogously, the *first excursion* $e_0 = \{e_0(n)\}_n := \{S_n\}_{0 \leq n \leq t_1}$ is distributed like the process $(S, \mathbb{P}_y^\searrow)$, that we introduced in Section 1.4, conditioned to have a fixed lifetime $\zeta = t_1$. It is therefore clear that, in order to extract the scaling limits of $\mathbb{P}_{N,y,\varepsilon}$ as $N \to \infty$, one has to combine a good control on the law of the zero level set $\mathcal{I}_N$ with the asymptotic properties of the excursions, and in this respect the usefulness of Theorem 1.3 emerges.

### 1.6. Outline of the paper

The exposition is organized as follows:

- in Section 2 we collect some preliminary facts;
- Section 3 contains the proof of Theorem 1.1 in the special case $x_N \equiv 0$;
- in Section 4 we complete the proof of Theorem 1.1, allowing for nonzero starting points;
- in Section 5 we study the random walk conditioned to die at zero and its counterpart for Lévy processes, proving Theorem 1.3;
- in Appendix A we present a path decomposition of the chain $(S, \mathbb{P}_y^+)$ at its overall minimum, together with the proof of some minor results;
- in Appendix B we prove that the function $V(\cdot)$ (resp. $V^\sim(\cdot)$) is invariant for the semigroup of the random walk killed when it first enters the negative (resp. nonpositive) half-line.

## 2. Some preliminary facts

Throughout the paper we use the notation $\alpha_n \sim \beta_n$ to indicate that $\alpha_n/\beta_n \to 1$ as $n \to \infty$. We recall that a positive sequence $d_n$ is said to be *regularly varying* of index $\alpha \in \mathbb{R}$ (this will be denoted by $d_n \in R_\alpha$) if $d_n \sim L_n n^\alpha$ as $n \to \infty$, where $L_n$ is *slowly varying* in that $L_{\lfloor tn \rfloor}/L_n \to 1$ as $n \to \infty$, for every $t > 0$. If $d_n$ is regularly varying with index $\alpha \neq 0$, up to asymptotic equivalence we will always assume that $d_n = d(n)$, with $d(\cdot)$ a continuous, strictly monotone function [4], Th. 1.5.3. Observe that if $d_n \in R_\alpha$ then $d^{-1}(n) \in R_{1/\alpha}$ and $1/d_n \in R_{-\alpha}$.

By the standard theory of stability, assumption (1.6) yields $a_n = a(n) \in R_{1/\alpha}$. In the following lemma we determine the asymptotic behaviour of the sequence $\mathbb{P}(\mathcal{C}_N)$ and of the function $V(\cdot)$, that will play a major role in the following sections.

**Lemma 2.1.** *The asymptotic behaviour of $V(x)$ and $\mathbb{P}(\mathcal{C}_N)$ are given by*

$$V(x) \sim C_1 \cdot c^{-1}(x), \quad x \to \infty, \qquad \mathbb{P}(\mathcal{C}_N) \sim \frac{C_2}{b^{-1}(N)}, \quad N \to \infty, \tag{2.1}$$

where $b(\cdot)$ and $c(\cdot)$ are continuous, strictly increasing functions such that $b(n) \in R_{1/\overline{\rho}}$ and $c(n) \in R_{1/\alpha\overline{\rho}}$. Moreover, $b(\cdot)$ and $c(\cdot)$ can be chosen such that $c = a \circ b$.

**Proof.** We recall that our random walk is attracted to a stable law of index $\alpha$ and positivity parameter $\rho$, and that we have set $\overline{\rho} := 1 - \rho$. Then by [15, 16, 23] we have that $\overline{T}_1$ and $\overline{H}_1$ are in the domain of attraction of the positive stable law of index respectively $\overline{\rho}$ and $\alpha\overline{\rho}$ (in the case $\alpha\overline{\rho} = 1$ by "the positive stable law of index 1" we simply mean the Dirac mass $\delta_1(\mathrm{d}x)$ at $x = 1$). The norming sequences of $\overline{T}$ and $\overline{H}$ will be denoted respectively by $b(n) \in R_{1/\overline{\rho}}$ and $c(n) \in R_{1/\alpha\overline{\rho}}$, where the functions $b(\cdot)$ and $c(\cdot)$ can be chosen continuous, increasing and such that $c = a \circ b$, cf. [15].

Recalling the definition (1.1), by standard Tauberian theorems (see [4], Section 8.2, for the $\alpha\overline{\rho} < 1$ case and [4], Section 8.8, for the $\alpha\overline{\rho} = 1$ case) we have that the asymptotic behaviour of $V(x)$ is given by

$$V(x) \sim C_1 \cdot c^{-1}(x), \quad x \to \infty, \tag{2.2}$$

where $C_1$ is a positive constant. In particular, $V(x) \in R_{\alpha\overline{\rho}}$.

Finally observe that since $\mathcal{C}_N = (\overline{T}_1 > N)$, the asymptotic behaviour of $\mathbb{P}(\mathcal{C}_N)$ is given by [20], Section XIII.6:

$$\mathbb{P}(\mathcal{C}_N) \sim C_2 / b^{-1}(N), \quad N \to \infty, \tag{2.3}$$

where $C_2$ is a positive constant. □

## 3. Convergence of $\mathbf{P}^{+,N}$

In this section we prove Theorem 1.1 in the special case when $x_N = x = 0$ for all $N$, i.e. we show that $\mathbf{P}^{+,N} \Rightarrow \mathbf{P}^+$ as $N \to \infty$.

### 3.1. Some preparation

We introduce the spaces $\Omega_t := D([0,t], \mathbb{R})$, $t > 0$ of càdlàg paths which are defined on the time interval $[0,t]$. For each $t$, the space $\Omega_t$ is endowed with the Skorohod topology, and with some misuse of notations we will call $\{\mathcal{F}_s\}_{s \in [0,t]}$ the natural filtration generated by the canonical process $X$ defined on this space.

We denote by $\mathbf{P}^{(m)}$ the law on $\Omega_1$ of the *meander* of length 1 associated to $(X, \mathbf{P})$, that is the rescaled post-minimum process of $(X, \mathbf{P})$, see [12]. It may also be defined more explicitly as the following weak limit:

$$\mathbf{P}^{(m)} = \lim_{x \downarrow 0} \mathbf{P}_x(\cdot | \underline{X}_1 \geq 0),$$

where $\underline{X}_1 = \inf_{0 \leq s \leq 1} X_s$, see Theorem 1 in [12]. Thus the law $\mathbf{P}^{(m)}$ may be considered as the law of the Lévy process $(X, \mathbf{P})$ conditioned to stay positive on the time interval $(0,1)$, whereas we have seen that the law $\mathbf{P}^+$ corresponds to an analogous conditioning but over the whole real half-line $(0, \infty)$. Actually it is proved in [12] that these measures are absolutely continuous with respect to each other: for every event $A \in \mathcal{F}_1$,

$$\mathbf{P}^+(A) = \mathbf{E}^{(m)}(U(X_1)\mathbf{1}_A), \tag{3.1}$$

where $U(x) := C_3 \cdot \widetilde{U}(x)$ and $C_3$ is a positive constant (the function $\widetilde{U}(x)$ has been defined in (1.4)).

Analogously we denote by $\mathbb{P}^{(m),N}$ the law on $\Omega_{\mathrm{RW}}$ corresponding to the random walk $(S, \mathbb{P})$ conditioned to stay nonnegative up to epoch $N$, that is,

$$\mathbb{P}^{(m),N} := \mathbb{P}(\cdot | \mathcal{C}_N).$$

As in the continuous setting, the two laws $\mathbb{P}^+$ and $\mathbb{P}^{(m),N}$ are mutually absolutely continuous: for every $B \in \sigma(S_1, \ldots, S_N)$ we have

$$\mathbb{P}^+(B) = \mathbb{P}(\mathcal{C}_N) \cdot \mathbb{E}^{(m),N}(V(S_N)\mathbf{1}_B), \tag{3.2}$$

where we recall that $V(x)$ defined in (1.1) is the renewal function of the strict descending ladder height process of the random walk. Note that in this case, relation (3.2) is a straightforward consequence of the definitions of the probability measures $\mathbb{P}^+$ and $\mathbb{P}^{(m),N}$.

Before getting into the proof we recall the invariance principle for the meander, that has been proven in more and more general settings in [6, 14] and [24]: introducing the rescaled meander measure $\mathbf{P}^{(m),N} := \mathbb{P}^{(m),N} \circ (\phi_N)^{-1}$ on $\Omega_1$ (here $\phi_N$ is to be understood as a map from $\Omega_{\mathrm{RW}}$ to $\Omega_1$), we have

$$\mathbf{P}^{(m),N} \implies \mathbf{P}^{(m)}, \quad N \to \infty. \tag{3.3}$$

*3.2. Proof of Theorem 1.1 ($x_N = x = 0$)*

Recalling the definition of $\mathbf{P}^{+,N}$ given in (1.8), from relation (3.2) we easily deduce the corresponding absolute continuity relation between $\mathbf{P}^{+,N}$ restricted to $\Omega_1$ and $\mathbf{P}^{(m),N}$: for every event $A \in \mathcal{F}_1$ we have

$$\mathbf{P}^{+,N}(A) = \mathbf{E}^{(m),N}(V_N(X_1)\mathbf{1}_A), \tag{3.4}$$

where we have introduced the rescaled renewal function

$$V_N(x) := \mathbb{P}(\mathcal{C}_N) \cdot V(a_N x). \tag{3.5}$$

We first prove that the sequence of measures $\mathbf{P}^{+,N}$ restricted to $\Omega_1$ converges weakly towards the measure $\mathbf{P}^+$ restricted to $\Omega_1$. To do so, we have to show that for every functional $H : \Omega_1 \to \mathbb{R}$ which is bounded and continuous one has $\mathbf{E}^{+,N}(H) \to \mathbf{E}^+(H)$ as $N \to \infty$. Looking at (3.1) and (3.4), this is equivalent to showing that

$$\mathbf{E}^{(m),N}(H \cdot V_N(X_1)) \to \mathbf{E}^{(m)}(H \cdot U(X_1)), \quad N \to \infty. \tag{3.6}$$

The basic idea is to show that $V_N(x) \to U(x)$ as $N \to \infty$ and then to use the invariance principle (3.3). However some care is needed, because the functions $V_N(\cdot)$ and $U(\cdot)$ are unbounded and the coordinate projections $X_t$ are not continuous in the Skorohod topology.

We start by introducing for $M > 0$ the cut function $I_M(x)$ (which can be viewed as a continuous version of $\mathbf{1}_{(-\infty,M]}(x)$):

$$I_M(x) := \begin{cases} 1, & x \leq M, \\ M+1-x, & M \leq x \leq M+1, \\ 0, & x \geq M+1. \end{cases} \tag{3.7}$$

The first step is to restrict the values of $X_1$ to a compact set. More precisely, we can decompose the l.h.s. of (3.6) as

$$\mathbf{E}^{(m),N}(H \cdot V_N(X_1)) = \mathbf{E}^{(m),N}(H \cdot V_N(X_1) \cdot I_M(X_1)) + \mathbf{E}^{(m),N}(H \cdot V_N(X_1) \cdot (1 - I_M(X_1))),$$

and analogously for the r.h.s. Then by the triangle inequality we easily get

$$|\mathbf{E}^{(m),N}(H \cdot V_N(X_1)) - \mathbf{E}^{(m)}(H \cdot U(X_1))|$$
$$\leq |\mathbf{E}^{(m),N}(H \cdot V_N(X_1) \cdot I_M(X_1)) - \mathbf{E}^{(m)}(H \cdot U(X_1) \cdot I_M(X_1))|$$
$$+ |\mathbf{E}^{(m),N}(H \cdot V_N(X_1) \cdot (1 - I_M(X_1)))| + |\mathbf{E}^{(m)}(H \cdot U(X_1) \cdot (1 - I_M(X_1)))|.$$

Since $H$ is bounded by some positive constant $C_4$ and the terms $V_N(X_1) \cdot (1 - I_M(X_1))$ and $U(X_1) \cdot (1 - I_M(X_1))$ are nonnegative, we get

$$|\mathbf{E}^{(m),N}(H \cdot V_N(X_1)) - \mathbf{E}^{(m)}(H \cdot U(X_1))|$$
$$\leq |\mathbf{E}^{(m),N}(H \cdot V_N(X_1) \cdot I_M(X_1)) - \mathbf{E}^{(m)}(H \cdot U(X_1) \cdot I_M(X_1))|$$
$$+ C_4 \mathbf{E}^{(m),N}(V_N(X_1) \cdot (1 - I_M(X_1))) + C_4 \mathbf{E}^{(m)}(U(X_1) \cdot (1 - I_M(X_1))).$$

However by definition we have $\mathbf{E}^{(m),N}(V_N(X_1)) = 1$ and $\mathbf{E}^{(m)}(U(X_1)) = 1$, hence

$$|\mathbf{E}^{(m),N}(H \cdot V_N(X_1)) - \mathbf{E}^{(m)}(H \cdot U(X_1))|$$
$$\leq |\mathbf{E}^{(m),N}(H \cdot V_N(X_1) \cdot I_M(X_1)) - \mathbf{E}^{(m)}(H \cdot U(X_1) \cdot I_M(X_1))|$$
$$+ C_4(1 - \mathbf{E}^{(m),N}(V_N(X_1) \cdot I_M(X_1))) + C_4(1 - \mathbf{E}^{(m)}(U(X_1) \cdot I_M(X_1))). \tag{3.8}$$

Next we claim that for every $M > 0$ the first term in the r.h.s. of (3.8) vanishes as $N \to \infty$, namely

$$|\mathbf{E}^{(m),N}(H \cdot V_N(X_1) \cdot I_M(X_1)) - \mathbf{E}^{(m)}(H \cdot U(X_1) \cdot I_M(X_1))| \to 0, \quad N \to \infty. \tag{3.9}$$

Observe that this equation yields also the convergence as $N \to \infty$ of the second term in the r.h.s. of (3.8) towards the third term (just take $H \equiv 1$), and note that the third term can be made arbitrarily small by choosing $M$ sufficiently large, again because $\mathbf{E}^{(m)}(U(X_1)) = 1$. Therefore from (3.9) it actually follows that the l.h.s. of (3.8) vanishes as $N \to \infty$, that is Eq. (3.6) holds true.

It remains to prove (3.9). By the triangle inequality we get

$$|\mathbf{E}^{(m),N}(H \cdot V_N(X_1) \cdot I_M(X_1)) - \mathbf{E}^{(m)}(H \cdot U(X_1) \cdot I_M(X_1))|$$
$$\leq C_4 \sup_{x \in [0,M]} |V_N(x) - U(x)| + |\mathbf{E}^{(m),N}(H \cdot U(X_1) \cdot I_M(X_1)) - \mathbf{E}^{(m)}(H \cdot U(X_1) \cdot I_M(X_1))|, \tag{3.10}$$

and we now show that both terms in the r.h.s. above vanish as $N \to \infty$.

By the uniform convergence property of regularly varying sequences [4], Theorem 1.2.1, it follows that for any $0 < \eta < M < \infty$

$$V(sx) = x^{\alpha \bar{\rho}} V(s)(1 + o(1)), \quad s \to \infty,$$

uniformly for $x \in [\eta, M]$ (recall that $V(s) \in R_{\alpha \bar{\rho}}$ as $s \to \infty$), hence from (3.5) we get

$$V_N(x) = (\mathbb{P}(\mathcal{C}_N) \cdot V(a_N)) x^{\alpha \bar{\rho}}(1 + o(1)), \quad N \to \infty, \tag{3.11}$$

uniformly for $x \in [\eta, M]$. Let us look at the prefactor above: by (2.2) we have

$$V(a_N) \sim C_1 \cdot c^{-1}(a_N) = C_1 \cdot b^{-1}(N), \quad N \to \infty,$$

where in the second equality we have used the fact that $c = a \circ b$ and hence $b^{-1} = c^{-1} \circ a$. Then it follows by (2.3) that $(\mathbb{P}(\mathcal{C}_N) \cdot V(a_N)) \to C_1 \cdot C_2$ as $N \to \infty$. In fact $C_1 \cdot C_2$ coincides with the constant $C_3$ defined after (3.1), hence we can rewrite (3.11) as

$$V_N(x) = U(x)(1 + o(1)), \quad N \to \infty,$$

uniformly for $x \in [\eta, M]$. However we are interested in absolute rather than relative errors, and using the fact that $V_N(\cdot)$ is increasing it is easy to throw away the $\eta$, getting that for every $M > 0$,

$$\sup_{x \in [0,M]} |V_N(x) - U(x)| \to 0, \quad N \to \infty. \tag{3.12}$$

Finally we are left with showing that the second term in the r.h.s. of (3.10) vanishes as $N \to \infty$. As already mentioned, the coordinate projections $X_t$ are not continuous in the Skorohod topology. However in our situation we have that $X_1 = X_{1-}$, $\mathbf{P}^+$-a.s., and this yields that the discontinuity set of the projection $X_1$ is $\mathbf{P}^+$-negligible. Therefore the functional

$$\Omega_1 \ni \omega \longmapsto H(\omega) \cdot U(X_1(\omega)) \cdot I_M(X_1(\omega))$$

is $\mathbf{P}^+$-a.s. continuous and bounded, and the conclusion follows directly from the invariance principle (3.3) for the meander.

Thus we have proved that the measure $\mathbf{P}^{+,N}$ restricted to $\Omega_1$ converges weakly to the measure $\mathbf{P}^+$ restricted to $\Omega_1$. Now it is not difficult to see that a very similar proof shows that $\mathbf{P}^{+,N}$ restricted to $\Omega_t$ converges weakly towards $\mathbf{P}^+$ restricted to $\Omega_t$, for each $t > 0$. Then it remains to apply Theorem 16.7 in [5] to obtain the weak convergence on the whole $\Omega$.

## 4. Convergence of $\mathbf{P}^{+,N}_{x_N}$

In this section we complete the proof of Theorem 1.1, allowing for nonzero starting points. We recall the laws defined in the introduction:

- $\mathbf{P}^+_x$ is the law on $\Omega$ of the Lévy process starting from $x \geq 0$ and conditioned to stay positive on $(0, \infty)$, cf. (1.5);
- $\mathbb{P}^+_y$ is the law on $\Omega_{\mathrm{RW}}$ of the random walk starting from $y \geq 0$ and conditioned to stay positive on $\mathbb{N}$, cf. (1.2);
- $\mathbf{P}^{+,N}_x$ is the law on $\Omega$ corresponding to the rescaling of $\mathbb{P}^+_y$, where $x = y/a_N$, cf. (1.8).

For ease of exposition we consider separately the cases $x > 0$ and $x = 0$.

*4.1. The case $x > 0$*

By the same arguments as in Section 3, we only need to prove the weak convergence of the sequence $\mathbf{P}^{+,N}_{x_N}$ restricted to $\Omega_1$ towards the measure $\mathbf{P}^+_x$ restricted to $\Omega_1$. Let $H : \Omega_1 \to \mathbb{R}$ be a continuous functional which is bounded by a constant $C$. Definitions (1.2), (1.5) and (1.8) give

$$|\mathbf{E}^{+,N}_{x_N}(H) - \mathbf{E}^+_x(H)| = |V_N(x_N)^{-1}\mathbf{E}^N_{x_N}(HV_N(X_1)\mathbf{1}_{\{\underline{X}_1 \geq 0\}}) - U(x)^{-1}\mathbf{E}_x(HU(X_1)\mathbf{1}_{\{\underline{X}_1 \geq 0\}})|.$$

Since Eq. (3.12) yields

$$V_N(x_N) \to U(x), \quad N \to \infty, \tag{4.1}$$

and since $U(x) > 0$ for $x > 0$, to obtain our result it suffices to show that

$$|\mathbf{E}^N_{x_N}(HV_N(X_1)\mathbf{1}_{\{\underline{X}_1 \geq 0\}}) - \mathbf{E}_x(HU(X_1)\mathbf{1}_{\{\underline{X}_1 \geq 0\}})| \to 0, \quad N \to \infty. \tag{4.2}$$

We proceed as in Section 3. It is easy to check that the discontinuity set of the functional $\mathbf{1}_{\{\underline{X}_1 \geq 0\}}$ is $\mathbf{P}_x$-negligible. Therefore the functional

$$\Omega_1 \ni \omega \longmapsto H(\omega) \cdot \mathbf{1}_{\{\underline{X}_1 \geq 0\}} \cdot U(X_1(\omega)) \cdot I_M(X_1(\omega))$$

is $\mathbf{P}^+$-a.s. continuous and bounded (we recall that the function $I_M(\cdot)$ has been defined in (3.7)). Hence, using the invariance principle of the unconditioned law, that is

$$\mathbf{P}^N_{x_N} \Longrightarrow \mathbf{P}_x, \quad N \to \infty, \tag{4.3}$$

see for instance [27], we deduce that for any $M > 0$ as $N \to \infty$

$$|\mathbf{E}^N_{x_N}(H \cdot U(X_1) \cdot \mathbf{1}_{\{\underline{X}_1 \geq 0\}} \cdot I_M(X_1)) - \mathbf{E}_x(H \cdot U(X_1) \cdot \mathbf{1}_{\{\underline{X}_1 \geq 0\}} \cdot I_M(X_1))| \to 0.$$

Then in complete analogy with (3.10), from the triangle inequality and (3.12) it follows that as $N \to \infty$

$$|\mathbf{E}^N_{x_N}(H \cdot V_N(X_1) \cdot \mathbf{1}_{\{\underline{X}_1 \geq 0\}} \cdot I_M(X_1)) - \mathbf{E}_x(H \cdot U(X_1) \cdot \mathbf{1}_{\{\underline{X}_1 \geq 0\}} \cdot I_M(X_1))| \to 0. \tag{4.4}$$

Now since $H$ is bounded by $C$, from the triangle inequality we can write

$$|\mathbf{E}^N_{x_N}(H \cdot V_N(X_1) \cdot \mathbf{1}_{\{\underline{X}_1 \geq 0\}}) - \mathbf{E}_x(H \cdot U(X_1) \cdot \mathbf{1}_{\{\underline{X}_1 \geq 0\}})|$$
$$\leq |\mathbf{E}^N_{x_N}(H \cdot V_N(X_1) \cdot \mathbf{1}_{\{\underline{X}_1 \geq 0\}} \cdot I_M(X_1)) - \mathbf{E}_x(H \cdot U(X_1) \cdot \mathbf{1}_{\{\underline{X}_1 \geq 0\}} \cdot I_M(X_1))|$$
$$+ C(V_N(x_N) - \mathbf{E}^N_{x_N}(V_N(X_1) \cdot \mathbf{1}_{\{\underline{X}_1 \geq 0\}} \cdot I_M(X_1))) + C(U(x) - \mathbf{E}_x(U(X_1) \cdot \mathbf{1}_{\{\underline{X}_1 \geq 0\}} \cdot I_M(X_1))), \tag{4.5}$$

where for the two last terms we have used the equalities

$$\mathbf{E}_{x_N}^N(V_N(X_1)\mathbf{1}_{\{\underline{X}_1 \geq 0\}}) = V_N(x_N) \quad \text{and} \quad \mathbf{E}_x(U(X_1)\mathbf{1}_{\{\underline{X}_1 \geq 0\}}) = U(x). \tag{4.6}$$

We deduce from (4.4) and (4.1) that when $N \to \infty$ the first term in the r.h.s. of (4.5) tends to 0 and that the second term converges towards the third one. Thanks to (4.6), the latter can be made arbitrarily small as $M \to \infty$, hence our result is proved in the case $x > 0$.

*4.2. The case $x = 0$*

We are going to follow arguments close to those developed by Bryn-Jones and Doney in the Gaussian case [7]. The proof uses the following path decomposition of $(X, \mathbf{P}_x^{+,N})$ at its overall minimum time, which is very similar to the analogous result for Lévy processes proven in [11]. The proofs of the next two lemmas are postponed to Appendix A.

**Lemma 4.1.** *Let $\mu = \inf\{t: X_t = \inf_{s \geq 0} X_s\}$. Then under $\mathbf{P}_x^{+,N}$, the post-minimum process $\{X_{t+\mu} - \inf_{s \geq 0} X_s, t \geq 0\}$ has law $\mathbf{P}_0^{+,N}$ and the overall minimum $\inf_{s \geq 0} X_s$ is distributed as*

$$\mathbf{P}_x^{+,N}\left(\inf_{s \geq 0} X_s \geq z\right) = \frac{V_N(x-z)}{V_N(x)}, \quad 0 \leq z \leq x. \tag{4.7}$$

*Moreover, under $\mathbf{P}_x^{+,N}$ the pre-minimum process $\{X_s, s \leq \mu\}$ and the post-minimum process are independent.*

**Lemma 4.2.** *Let $\mu = \inf\{t: X_t = \inf_{s \geq 0} X_s\}$. If $x_N \to 0$ as $N \to \infty$, then under $\mathbf{P}_{x_N}^{+,N}$ the maximum of the process before time $\mu$ and the time $\mu$ itself converge in probability to 0, i.e. for all $\varepsilon > 0$*

$$\lim_{N \to +\infty} \mathbf{P}_{x_N}^{+,N}(\mu \geq \varepsilon) = 0 \quad \text{and} \quad \lim_{N \to +\infty} \mathbf{P}_{x_N}^{+,N}\left(\sup_{0 \leq s \leq \mu} X_s \geq \varepsilon\right) = 0. \tag{4.8}$$

Now let $(\Omega', \mathcal{F}', P)$ be a probability space on which are defined the processes $X^{(N)}$, $N \in \mathbb{N} \cup \{\infty\}$, such that $X^{(N)}$ has law $\mathbf{P}_0^{+,N}$, $X^{(\infty)}$ has law $\mathbf{P}_0^+$ and

$$X^{(N)} \to X^{(\infty)} \quad P\text{-almost surely}. \tag{4.9}$$

Since we have already proved that $\mathbf{P}_0^{+,N} \Rightarrow \mathbf{P}_0^+$ as $N \to \infty$, such a construction is possible in virtue of Skorohod's Representation Theorem. We assume that on the same space is defined a sequence of processes $Y^{(N)}$, $N \in \mathbb{N}$, such that $Y^{(N)}$ has law $\mathbf{P}_{x_N}^{+,N}$ and is independent of $X^{(N)}$. Then we set $\mu^{(N)} = \inf\{t: Y_t^{(N)} = \inf_{s \geq 0} Y_s^{(N)}\}$ and for each $N$ we define

$$Z_t^{(N)} := \begin{cases} Y_t^{(N)} & \text{if } t < \mu^{(N)}, \\ X_{t-\mu^{(N)}}^{(N)} & \text{if } t \geq \mu^{(N)}. \end{cases}$$

It follows from Lemma 4.1 that $Z^{(N)}$ has law $\mathbf{P}_{x_N}^{+,N}$. Moreover, from (4.8) we deduce that the process $\{Y_t^{(N)}\mathbf{1}_{\{t < \mu^{(N)}\}}, t \geq 0\}$ converges in $P$-probability as $N \to \infty$ towards the process which is identically equal to 0 in the Skorohod's space. Combining this result with the almost sure convergence (4.9), we deduce that for every subsequence $(N_k)$ there exists a sub-subsequence $(N'_k)$ such that $P$-a.s. $Z^{(N'_k)} \to X^{(\infty)}$ as $k \to \infty$ in the Skorohod's topology, and this completes the proof.

## 5. Conditioning to die at zero

In this section we study in detail the conditioning to die at zero for random walks and Lévy processes, that has been already introduced in Section 1.4. We will be dealing with Markov processes with values in $\mathbb{R}^+$ and having 0 as absorbing state. With some abuse of notation, the first hitting time of 0 by a path in $\Omega$ or in $\Omega_{\mathrm{RW}}$ will be denoted in both cases by $\zeta$, that is:

$$\zeta := \inf\{n \in \mathbb{Z}^+ \colon S_n = 0\} \quad \text{and} \quad \zeta := \inf\{t \geq 0 \colon X_t = 0\}.$$

For any interval $I$ of $\mathbb{R}$, we denote the first hitting time of $I$ by the canonical processes $S$ and $X$ in $\Omega$ and $\Omega_{\mathrm{RW}}$ respectively by

$$T_I = \inf\{n \in \mathbb{Z}^+ \colon S_n \in I\} \quad \text{and} \quad \tau_I = \inf\{t \geq 0 \colon X_t \in I\}.$$

### 5.1. Lévy process conditioned to die at zero

We now introduce the conditioning to die at zero for Lévy processes that has been studied in [11], Section 4. We recall from Section 1.2 that the function $\widetilde{U}(x) = x^{\alpha\overline{\rho}}$ is invariant for the semigroup of $(X, \mathbf{P}_x)$ killed at time $\tau_{(-\infty,0]}$. (Note that for Lévy processes $\tau_{(-\infty,0]} = \tau_{(-\infty,0)}$, a.s.) In the setting of this section, the killed process is identically equal to 0 after time $\tau_{(-\infty,0]}$. The derivative $\widetilde{U}'(x) = \alpha\overline{\rho} x^{\alpha\overline{\rho}-1}$ is excessive for the same semigroup, that is for every $t \geq 0$

$$\widetilde{U}'(x) \geq \mathbf{E}_x(\widetilde{U}'(X_t)\mathbf{1}_{(\tau_{(-\infty,0]}>t)}).$$

Then one can define the $h$-transform of the killed process by the function $\widetilde{U}'(\cdot)$, that is the Markovian law $\mathbf{P}_x^{\searrow}$ on $\Omega$ defined for $x > 0$, $t > 0$ and $A \in \mathcal{F}_t$ by

$$\mathbf{P}_x^{\searrow}(A, \zeta > t) := \frac{1}{\widetilde{U}'(x)} \mathbf{E}_x(\widetilde{U}'(X_t)\mathbf{1}_A \mathbf{1}_{(\tau_{(-\infty,0]}>t)}). \tag{5.1}$$

Note that from [13], XVI.30, the definition (5.1) is still valid when replacing $t$ by any stopping time of the filtration $(\mathcal{F}_t)$. Since 0 is an absorbing state, Eq. (5.1) entirely determines the law $\mathbf{P}_x^{\searrow}$, in particular $\mathbf{P}_0^{\searrow}$, is the law of the degenerated process $X \equiv 0$. The process $(X, \mathbf{P}_x^{\searrow})$ is called the Lévy process conditioned to die at zero. This terminology is justified by the following result, proven in [11], Proposition 3: for all $x, \beta, t > 0$ and $A \in \mathcal{F}_t$ one has

$$\lim_{\varepsilon \to 0} \mathbf{P}_x(A, \tau_{(-\infty,\beta)} > t | \underline{X}_{\tau_{(-\infty,0)-}} \leq \varepsilon) = \mathbf{P}_x^{\searrow}(A, \tau_{[0,\beta)} > t). \tag{5.2}$$

We also emphasize that the process $(X, \mathbf{P}_x^{\searrow})$ a.s. hits 0 in a finite time and that either it has a.s. no negative jumps or it reaches 0 by an accumulation of negative jumps, i.e. $\mathbf{P}_x^{\searrow}(\zeta < \infty, X_{\zeta-} = 0) = 1$.

### 5.2. Random walk conditioned to die at zero

We recall the construction in the random walk setting, that we started in Section 1.4. We assume that $S$ is $\mathbb{Z}$-valued and aperiodic and we introduce the Markovian family of laws $\mathbb{P}_y^{\searrow}$, $y \in \mathbb{R}_+$, on $\Omega_{\mathrm{RW}}$ defined for $N \in \mathbb{N}$ and for $B \in \sigma(S_1, \ldots, S_N)$ by

$$\mathbb{P}_y^{\searrow}(B, \zeta > N) := \mathbb{P}_y(B, T_{(-\infty,0]} > N | S_{T_{(-\infty,0]}} \in (-1,0]), \quad y > 0, \tag{5.3}$$

while $\mathbb{P}_0^{\searrow}$ is the law of the process $S \equiv 0$. We point out that the hypothesis of aperiodicity ensures that for all $y > 0$

$$W(y) := \mathbb{P}_y(S_{T_{(-\infty,0]}} \in (-1,0]) > 0, \tag{5.4}$$

so that the conditioning in (5.3) makes sense. To prove this relation, first notice that $W(y) = W(\lceil y \rceil)$, where $\lceil y \rceil$ denotes the upper integer part of $y$. Moreover, for $n \in \mathbb{N}$ an inclusion lower bound together with the Markov property yields

$$W(n) = \mathbb{P}_n(S_{T_{(-\infty,0]}} = 0) \geq \mathbb{P}_n\left(\bigcap_{i=0}^{n-1} S_{T_{(-\infty,i]}} = i\right) = (\mathbb{P}_1(S_{T_{(-\infty,0]}} = 0))^n = (W(1))^n,$$

hence we are left with showing that $W(1) = \mathbb{P}_0(\overline{H}_1 = 1) > 0$ (we recall that $\overline{H}_1$ is the first descending ladder height, defined in Section 1.1). To this purpose, we use a basic combinatorial identity for general random walks discovered by Alili and Doney, cf. Eq. (6) in [1], that in our case gives

$$\mathbb{P}_0(\overline{H}_1 = 1, \overline{T}_1 = n) = \frac{1}{n}\mathbb{P}_0(S_n = -1).$$

It only remains to observe that Gnedenko's Local Limit Theorem [4], Theorem 8.4.1, yields the positivity of the r.h.s. for large $n$. One can easily check that for $y > 0$

$$W(y) = \mathbb{P}_0(S_{T_{(-\infty,-y]}} \in (-y-1,-y]) = \mathbb{P}_0(\exists k \colon \overline{H}_k = \lceil y \rceil) = V(\lceil y \rceil) - V(\lceil y \rceil - 1), \tag{5.5}$$

where $V(\cdot)$ is the renewal function of the renewal process $\{\overline{H}_k\}$, as defined in (1.1).

The following lemma gives a useful description of $(S, \mathbb{P}_y^{\searrow})$ as an $h$-transform.

**Lemma 5.1.** *The Markov chain $(S, \mathbb{P}_y^{\searrow})$, $y \geq 0$, is an $h$-transform of $(S, \mathbb{P}_y)$ killed when it enters the nonpositive half-line $(-\infty, 0]$ corresponding to the excessive function $W(y)$, $y \geq 0$, i.e. for any $N \in \mathbb{N}$ and for any $B \in \sigma(S_1, \ldots, S_N)$*

$$\mathbb{P}_y^{\searrow}(B, \zeta > N) = \frac{1}{W(y)} \mathbb{E}_y(B \mathbf{1}_{\{T_{(-\infty,0]} > N\}} W(S_N)). \tag{5.6}$$

**Proof.** It is just a matter of applying the definition (5.3) and the Markov property, getting for $N \in \mathbb{N}$ and for $B \in \sigma(S_1, \ldots, S_N)$

$$\mathbb{P}_y^{\searrow}(B, \zeta > N) = \frac{\mathbb{E}_y(\mathbf{1}_B \mathbf{1}_{\{T_{(-\infty,0]} > N\}} \mathbb{P}_{S_N}(S_{T_{(-\infty,0]}} \in (-1,0]))}{\mathbb{P}_y(S_{T_{(-\infty,0]}} \in (-1,0])} = \frac{1}{W(y)} \mathbb{E}_y(\mathbf{1}_B \mathbf{1}_{\{T_{(-\infty,0]} > N\}} W(S_N)), \tag{5.7}$$

which also shows that the function $W(\cdot)$ is indeed excessive for $(S, \mathbb{P}_y)$ killed when it enters the nonpositive half-line. $\square$

We point out that the special choice of the law of $(S, \mathbb{P})$ in $\mathbb{Z}$ has been done only in the aim of working in a simpler setting. However, it is clear from our construction that the conditioning to die at 0 may be defined for general laws with very few assumptions.

### 5.3. Proof of Theorem 1.3

We start observing that by the definition of $\mathbf{P}_x^{\searrow}$, see (1.10), for any $t > 0$ and any $\mathcal{F}_t$-measurable functional $F$ one has

$$\mathbf{E}_{x_N}^{\searrow, N}(F \mathbf{1}_{(\zeta > t)}) = \frac{1}{W_N(x_N)} \mathbf{E}_{x_N}^N(F \mathbf{1}_{(\tau_{(-\infty,0]} > t)} W_N(X_t)), \tag{5.8}$$

where we have introduced the rescaled function

$$W_N(x) := \frac{W(a_N x)}{W(a_N)}. \tag{5.9}$$

To lighten the exposition, we will limit ourselves to the case $\alpha\overline{\rho} < 1$ (the case $\alpha\overline{\rho} = 1$ is analogous but has to be handled separately). We stress that we are working under the additional hypothesis

$$W(x) \sim \frac{\alpha\overline{\rho}}{L(x)} x^{\alpha\overline{\rho}-1}, \quad x \to \infty, \tag{5.10}$$

see Section 1.4 for more details. Then by the Uniform Convergence Theorem for regularly varying functions with negative index we have that for every $\delta > 0$

$$\sup_{z \in [\delta, +\infty)} |W_N(z) - z^{\alpha\overline{\rho}-1}| \to 0, \quad N \to \infty, \tag{5.11}$$

cf. [4], Theorem 1.5.2. Moreover, we introduce the functions $\overline{W}(z) := \sup_{y \in [z,\infty)} W(y)$ and $\underline{W}(z) := \inf_{y \in [0,z]} W(y)$, and we have the following relations

$$W(z) \sim \overline{W}(z), \qquad W(z) \sim \underline{W}(z), \quad z \to \infty, \tag{5.12}$$

which follow from [4], Th. 1.5.3.

For ease of exposition we divide the rest of the proof in two parts, considering separately the cases $x > 0$ and $x = 0$.

*The case $x > 0$*

We will first show that for every $u, v$ such that $0 < u < v < x$ and for every bounded, continuous and $\mathcal{F}_\infty$-measurable functional $H$ one has

$$\mathbf{E}_{x_N}^{\searrow,N}(H(X^{(u,v)})\mathbf{1}_{(\zeta > \tau_{[u,v]})}) \to \mathbf{E}_x^{\searrow}(H(X^{(u,v)})\mathbf{1}_{(\zeta > \tau_{[u,v]})}), \quad N \to \infty, \tag{5.13}$$

where $X^{(u,v)} = (X_t \mathbf{1}_{\{t \leq \tau_{[u,v]}\}}, t \geq 0)$. Note that $H(X^{(u,v)})$ is $\mathcal{F}_{\tau_{[u,v]}}$-measurable. Moreover, since (5.6) is still valid when replacing $N$ by any stopping time of the filtration $(\sigma(S_1, \ldots, S_k))_k$, one easily checks that (5.8) extends to the first passage time $\tau_{[u,v]}$. Since $W_N(x_N) \to x^{\alpha\overline{\rho}-1} > 0$ by (5.11), it suffices to show that

$$\mathbf{E}_{x_N}^N(H(X^{(u,v)})W_N(X_{\tau_{[u,v]}})\mathbf{1}_{(\tau_{(-\infty,0]} > \tau_{[u,v]})}) \to \mathbf{E}_x(H(X^{(u,v)})(X_{\tau_{[u,v]}})^{\alpha\overline{\rho}-1}\mathbf{1}_{(\tau_{(-\infty,0]} > \tau_{[u,v]})}),$$

as $N \to \infty$. By the triangle inequality we obtain

$$|\mathbf{E}_{x_N}^N(H(X^{(u,v)})W_N(X_{\tau_{[u,v]}})\mathbf{1}_{(\tau_{(-\infty,0]} > \tau_{[u,v]})}) - \mathbf{E}_x(H(X^{(u,v)})(X_{\tau_{[u,v]}})^{\alpha\overline{\rho}-1}\mathbf{1}_{(\tau_{(-\infty,0]} > \tau_{[u,v]})})|$$
$$\leq |\mathbf{E}_{x_N}^N(H(X^{(u,v)})\mathbf{1}_{(\tau_{(-\infty,0]} > \tau_{[u,v]})}|W_N(X_{\tau_{[u,v]}}) - (X_{\tau_{[u,v]}})^{\alpha\overline{\rho}-1}|)|$$
$$+ |\mathbf{E}_{x_N}^N(H(X^{(u,v)})(X_{\tau_{[u,v]}})^{\alpha\overline{\rho}-1}\mathbf{1}_{(\tau_{(-\infty,0]} > \tau_{[u,v]})})$$
$$- \mathbf{E}_x(H(X^{(u,v)})(X_{\tau_{[u,v]}})^{\alpha\overline{\rho}-1}\mathbf{1}_{(\tau_{(-\infty,0]} > \tau_{[u,v]})})|. \tag{5.14}$$

By the right continuity of the canonical process, we have $X_{\tau_{[u,v]}} \in [u, v]$ a.s. on the event $(\tau_{(-\infty,0]} > \tau_{[u,v]})$. Moreover, since the functional $H$ is bounded by the positive constant $C_1$ we can estimate the first term in the r.h.s. above by

$$|\mathbf{E}_{x_N}^N(H(X^{(u,v)})\mathbf{1}_{(\tau_{(-\infty,0]} > \tau_{[u,v]})}|W_N(X_{\tau_{[u,v]}}) - (X_{\tau_{[u,v]}})^{\alpha\overline{\rho}-1}|)| \leq C_1 \sup_{z \in [u,v]} |W_N(z) - z^{\alpha\overline{\rho}-1}|,$$

which vanishes as $N \to \infty$ by Eq. (5.11). Moreover the second term in the r.h.s. of (5.14) tends to 0 as $N \to \infty$ thanks to the invariance principle $\mathbf{P}_{x_N}^N \Rightarrow \mathbf{P}_x$ for the unconditioned process, because the functional $(X_{\tau_{[u,v]}})^{\alpha\overline{\rho}-1}\mathbf{1}_{(\zeta > \tau_{[u,v]})}$ is bounded and its discontinuity set has zero $\mathbf{P}_x$-probability. This completes the proof of (5.13).

To obtain our result, we have to show that the left member of the inequality

$$|\mathbf{E}_{x_N}^{\searrow,N}(H) - \mathbf{E}_x^{\searrow}(H)|$$
$$\leq |\mathbf{E}_{x_N}^{\searrow,N}(H\mathbf{1}_{(\zeta > \tau_{[u,v]})}) - \mathbf{E}_x^{\searrow}(H\mathbf{1}_{(\zeta > \tau_{[u,v]})})| + |\mathbf{E}_{x_N}^{\searrow,N}(H\mathbf{1}_{(\zeta \leq \tau_{[u,v]})}) - \mathbf{E}_x^{\searrow}(H\mathbf{1}_{(\zeta \leq \tau_{[u,v]})})|, \tag{5.15}$$

tends to 0 as $N \to +\infty$. Fix $\varepsilon > 0$. Proposition 2 of [11] ensures that $\mathbf{P}_x^{\searrow}(\zeta > \tau_{(0,y]}) = 1$ for all $y > 0$, hence for all $v \in (0, x)$, there exists $u \in (0, v)$ such that $\mathbf{P}_x^{\searrow}(\zeta > \tau_{[u,v]}) \geq 1 - \varepsilon$. Moreover, from (5.13), $\mathbf{P}_{x_N}^{\searrow,N}(\zeta >$ $\tau_{[u,v]}) \to \mathbf{P}_x^{\searrow}(\boldsymbol{\zeta} > \tau_{[u,v]}))$, hence there exists $N_0$, such that for any $N \geq N_0$, $\mathbf{P}_{x_N}^{\searrow,N}(\boldsymbol{\zeta} > \tau_{[u,v]}) \geq 1 - 2\varepsilon$. So we have proved that for all $v \in (0, x)$, there exist $u \in (0, v)$ and $N_0$ such that $\forall N \geq N_0$

$$|\mathbf{E}_{x_N}^{\searrow,N}(H\mathbf{1}_{(\boldsymbol{\zeta} \leq \tau_{[u,v]})}) - \mathbf{E}_x^{\searrow}(H\mathbf{1}_{(\boldsymbol{\zeta} \leq \tau_{[u,v]})})| \leq 3C_1\varepsilon. \tag{5.16}$$

Now to deal with the first term of inequality (5.15), write

$$|\mathbf{E}_{x_N}^{\searrow,N}(H\mathbf{1}_{(\boldsymbol{\zeta} > \tau_{[u,v]})}) - \mathbf{E}_x^{\searrow}(H\mathbf{1}_{(\boldsymbol{\zeta} > \tau_{[u,v]})})|$$
$$\leq |\mathbf{E}_{x_N}^{\searrow,N}(H(X^{(u,v)})\mathbf{1}_{(\boldsymbol{\zeta} > \tau_{[u,v]})}) - \mathbf{E}_x^{\searrow}(H(X^{(u,v)})\mathbf{1}_{(\boldsymbol{\zeta} > \tau_{[u,v]})})|$$
$$+ \mathbf{E}_{x_N}^{\searrow,N}(|H - H(X^{(u,v)})|\mathbf{1}_{(\boldsymbol{\zeta} > \tau_{[u,v]})}) + \mathbf{E}_x^{\searrow}(|H - H(X^{(u,v)})|\mathbf{1}_{(\boldsymbol{\zeta} > \tau_{[u,v]})}). \tag{5.17}$$

The first term of the r.h.s. of (5.17) tends to 0 as $N \to \infty$, as we already proved above. It remains to analyze the other two terms. To this aim, let us consider on the space $D([0,\infty))$ a distance $d(\cdot, \cdot)$ that induces the Skorohod topology, e.g. as defined in [19], Section 3.5. We can choose it such that for $\xi, \eta \in D([0,\infty))$ we have $d(\xi, \eta) \leq \|\xi - \eta\|_\infty$, where $\|\cdot\|_\infty$ denotes the supremum norm over the real half-line. We can assume moreover that $H$ is a Lipschitz functional on $D([0,\infty))$:

$$|H(\xi)| \leq C_1, \qquad |H(\xi) - H(\eta)| \leq C_2 d(\xi, \eta) \quad \forall \xi, \eta \in D([0,\infty)), \tag{5.18}$$

because Lipschitz functionals determine convergence in law. Then setting $A_\varepsilon^{u,v} := (\sup_{t \in [\tau_{[u,v]}, \boldsymbol{\zeta})} X_t \leq \varepsilon)$, where $\varepsilon$ has been fixed above, and using (5.18) we have

$$\mathbf{E}_{x_N}^{\searrow,N}(|H - H(X^{(u,v)})|\mathbf{1}_{(\boldsymbol{\zeta} > \tau_{[u,v]})})$$
$$= \mathbf{E}_{x_N}^{\searrow,N}(|H - H(X^{(u,v)})|\mathbf{1}_{(\boldsymbol{\zeta} > \tau_{[u,v]})}\mathbf{1}_{A_\varepsilon^{u,v}}) + \mathbf{E}_{x_N}^{\searrow,N}(|H - H(X^{(u,v)})|\mathbf{1}_{(\boldsymbol{\zeta} > \tau_{[u,v]})}\mathbf{1}_{(A_\varepsilon^{u,v})^c})$$
$$\leq C_2\varepsilon \mathbf{P}_{x_N}^{\searrow,N}(\boldsymbol{\zeta} > \tau_{[u,v]}, A_\varepsilon^{u,v}) + 2C_1 \mathbf{P}_{x_N}^{\searrow,N}(\boldsymbol{\zeta} > \tau_{[u,v]}, (A_\varepsilon^{u,v})^c),$$

where we have used that $|H(X^{(u,v)}) - H| \leq C_2\|X^{(u,v)} - X\|_\infty = C_2 \sup_{t \in [\tau_{[u,v]}, \boldsymbol{\zeta})} |X_t|$ together with the fact that on the event $(\boldsymbol{\zeta} > \tau_{[u,v]}, A_\varepsilon^{u,v})$ one has by construction $\|X^{(u,v)} - X\|_\infty \leq \varepsilon$. An analogous expression can be derived under $\mathbf{P}_x^{\searrow}$.

To complete the proof it remains to show that one can find $u$ and $v$ sufficiently small such that

$$\mathbf{P}_x^{\searrow}(\boldsymbol{\zeta} > \tau_{[u,v]}, (A_\varepsilon^{u,v})^c) < \varepsilon \quad \text{and} \quad \limsup_N \mathbf{P}_{x_N}^{\searrow,N}(\boldsymbol{\zeta} > \tau_{[u,v]}, (A_\varepsilon^{u,v})^c) < \varepsilon. \tag{5.19}$$

Using the strong Markov property of $(X, \mathbf{P}_x^{\searrow})$ at time $\tau_{[u,v]}$, we obtain

$$\mathbf{P}_x^{\searrow}(\boldsymbol{\zeta} > \tau_{[u,v]}, (A_\varepsilon^{u,v})^c) = \mathbf{E}_x^{\searrow}\left(\mathbf{1}_{(\boldsymbol{\zeta} > \tau_{[u,v]})} \mathbf{P}_{X_{\tau_{[u,v]}}}^{\searrow}\left(\sup_{t \in [0,\boldsymbol{\zeta})} X_t > \varepsilon\right)\right). \tag{5.20}$$

But we easily check using (5.1) at the time $\tau_{(\varepsilon, \infty)}$ that

$$\mathbf{P}_z^{\searrow}\left(\sup_{t \in [0,\boldsymbol{\zeta})} X_t > \varepsilon\right) = z^{1 - \alpha\overline{\rho}} \mathbf{E}_z(X_{\tau_{(\varepsilon,\infty)}}^{\alpha\overline{\rho} - 1} \mathbf{1}_{(\tau_{(\varepsilon,\infty)} < \tau_{(-\infty,0]})}) \leq z^{1-\alpha\overline{\rho}} \varepsilon^{\alpha\overline{\rho} - 1} \to 0$$

as $z \to 0$. Coming back to (5.20), since $X_{\tau_{[u,v]}} \leq v$ a.s. on $(\boldsymbol{\zeta} > \tau_{[u,v]})$, we can find $v$ sufficiently small such that the first inequality in (5.19) holds for all $u \in (0, v)$. To obtain the second inequality, we use the Markov property of $(S, \mathbb{P}_{a_N x_N}^{\searrow})$ at time $T_{[a_N u, a_N v]}$ and Eq. (5.6) at time $T_{(\varepsilon a_N, \infty)}$, to obtain for all $u < v$

$$\mathbf{P}_{x_N}^{\searrow,N}(\boldsymbol{\zeta} > \tau_{[u,v]}, (A_\varepsilon^{u,v})^c) = \mathbb{E}_{a_N x_N}^{\searrow}\left(\mathbf{1}_{(T_{(-\infty,0]} > T_{[u a_N, v a_N]})} \frac{\mathbb{E}_{S_{T_{[u a_N, v a_N]}}}(\mathbf{1}_{(T_{(\varepsilon a_N, \infty)} < T_{(-\infty,0]})}, W(S_{T_{(\varepsilon a_N, \infty)}}))}{W(S_{T_{[u a_N, v a_N]}})}\right)$$
$$\leq \frac{\overline{W}(\varepsilon a_N)}{\underline{W}(v a_N)} \to \frac{v^{1-\alpha\overline{\rho}}}{\varepsilon^{1-\alpha\overline{\rho}}}, \quad N \to \infty,$$

where we have used (5.12) and the fact that $W(\cdot) \in R_{\alpha\overline{\rho}-1}$. Then it suffices to choose $v = \varepsilon^{1+1/(1-\alpha\overline{\rho})}$ to get the second inequality in (5.19) for all $u \in (0,v)$.

*The case $x = 0$*
Since the measure $\mathbf{P}_0^{\searrow}$ is the law of the process which is identically equal to zero, we only need to show that the overall maximum of the process $\mathbf{P}_{x_N}^{\searrow,N}$ converges to zero in probability, that is for every $\varepsilon > 0$

$$\lim_{N \to \infty} \mathbf{P}_{x_N}^{\searrow,N}\left(\sup_{t \in [0,\zeta)} X_t > \varepsilon\right) = 0. \tag{5.21}$$

It is convenient to rephrase this statement in terms of the random walk:

$$\mathbf{P}_{x_N}^{\searrow,N}\left(\sup_{t \in [0,\zeta)} X_t > \varepsilon\right) = \mathbb{P}_{x_N a_N}^{\searrow}\left(\sup_{n \in \{0,\zeta-1\}} S_n > \varepsilon a_N\right) = \mathbb{P}_{x_N a_N}^{\searrow}(T_{(\varepsilon a_N, \infty)} < \zeta).$$

From the definition of $\mathbb{P}_y^{\searrow}$ we can write

$$\mathbb{P}_{x_N a_N}^{\searrow}(T_{(\varepsilon a_N, \infty)} < \zeta) = \frac{1}{W(x_N a_N)} \mathbb{E}_{x_N a_N}\left(\mathbf{1}_{(T_{(\varepsilon a_N, \infty)} < T_{(-\infty, 0]})} W(S_{T_{(\varepsilon a_N, \infty)}}) \right) \le \frac{\overline{W}(\varepsilon a_N)}{W(x_N a_N)},$$

because by definition $S_{T_{(\varepsilon a_N, \infty)}} \ge \varepsilon a_N$ (we recall that $\overline{W}(z) := \sup_{y \ge z} W(y)$). Now it remains to show that

$$\lim_{N \to \infty} \frac{\overline{W}(\varepsilon a_N)}{W(x_N a_N)} \to 0. \tag{5.22}$$

Note that $x_N \to 0$ while $a_N \to \infty$, so that the asymptotic behaviour of $(x_N a_N)_N$ is not determined a priori. However, if we can show that (5.22) holds true whenever the sequence $(x_N a_N)_N$ has a (finite or infinite) limit, then the result in the general case will follow by a standard subsequence argument. Therefore we assume that $x_N a_N \to \ell \in [0, \infty]$ as $N \to \infty$. If $\ell < \infty$ then $\liminf_N W(x_N a_N) > 0$, and since $\overline{W}(\varepsilon a_N) \to 0$ Eq. (5.22) follows. On the other hand, if $\ell = \infty$ we have $x_N a_N \to \infty$ and to prove (5.22), from (5.12), we can replace $W(\cdot)$ by $\overline{W}(\cdot)$ which has the advantage of being decreasing. Since $x_N \to 0$, for every fixed $\delta > 0$ we have $x_N \le \delta$ for large $N$, hence from the monotonicity of $\overline{W}(\cdot)$ we get

$$\limsup_{N \to \infty} \frac{\overline{W}(a_N)}{\overline{W}(x_N a_N)} \le \limsup_{N \to \infty} \frac{\overline{W}(a_N)}{\overline{W}(\delta a_N)} = \delta^{1-\alpha\overline{\rho}},$$

where the last equality is nothing but the characteristic property of regularly varying functions. Since $\delta$ can be taken arbitrarily small, Eq. (5.22) is proven.

***Remark 5.2.*** *It follows from [11], Theorem 4 and Nagasawa's theory of time reversal that the returned process $(X_{(\zeta-t)-}, 0 \le t < \zeta)$ under $\mathbf{P}_x^{\searrow}$ has the same law as an h-transform of $(X, \mathbf{P}^{*,+})$, where $\mathbf{P}^{*,+}$ is the law of the process $(-X, \mathbf{P})$ conditioned to stay positive as defined in Section 1.2. Roughly speaking, it corresponds to $(X, \mathbf{P}^{*,+})$ conditioned to end at $x$. More specifically, if $p_t^{*,+}(y,z)$ stands for the semigroup of $(X, \mathbf{P}^{*,+})$, then it is the Markov process issued from 0 and with semigroup*

$$p_t^h(y,z) = \frac{h(z)}{h(y)} p_t^{*,+}(y,z),$$

*where $h(z) = \int_0^\infty p_t^{*,+}(z,x)\,\mathrm{d}t$. The same relationship between $(S, \mathbb{P}_y^{\searrow})$ and $(S, \mathbb{P}^{*,+})$ (where $\mathbb{P}^{*,+}$ is the law of the process $(-S, \mathbb{P})$ conditioned to stay positive) and the invariance principle established in Section 3, may provide another mean to obtain the main result of this section. The problem in this situation would reduce to the convergence of the discrete time equivalent of the function h.*

## Appendix A. Decomposition at the minimum for $\mathbb{P}_x^+$

*A.1. Proof of Lemma 4.1*

We start by rephrasing the Lemma in the space $\Omega_{\mathrm{RW}}$. We only assume that $(S,\mathbb{P})$ is a random walk that does not drift to $-\infty$, i.e. $\mathbb{P}(\limsup_{n\to\infty} S_n = +\infty) = 1$, and not that it is in the domain of attraction of a stable law. We denote by $(S, \mathbb{P}_y^+)$ the random walk started at $y \geq 0$ and conditioned to stay positive, defined in Section 1.1. Let $m = \inf\{n: S_n = \inf_{k \geq 0} S_k\}$ be the first time at which $S$ reaches its overall minimum. We are going to prove that under $\mathbb{P}_y^+$ the post-minimum process $\{S_{m+k} - S_m, k \geq 0\}$ has law $\mathbb{P}_0^+$ and is independent of the pre-minimum process $\{S_k, k \leq m\}$, and that the distribution of $S_m = \inf_{k \geq 0} S_k$ is given by

$$\mathbb{P}_y^+\left(\inf_{k\geq 0} S_k \geq x\right) = \frac{V(y-x)}{V(y)}, \quad 0 \leq x \leq y. \tag{A.1}$$

We start by proving the following basic relation: for every $A \in \sigma(S_n, n \in \mathbb{N})$ and for every $y \geq 0$

$$\mathbb{P}_y^+\left(A+y, \inf_{k\geq 0} S_k = y\right) = \frac{1}{V(y)} \mathbb{P}_0^+(A), \tag{A.2}$$

where the event $A+y$ is defined by $(S \in A+y) := (S - y \in A)$. By a standard monotone class argument, it suffices to prove this relation when $A \in \sigma(S_1, \ldots, S_N)$ for an arbitrary $N \in \mathbb{N}$. Then by the definition (1.2) of $\mathbb{P}_y^+$ for $n \geq N$ we can write

$$\mathbb{P}_y^+(A+y, S_1 \geq y, \ldots, S_n \geq y) = \frac{1}{V(y)} \mathbb{E}_y(V(S_n) \mathbf{1}_{(A+y, S_1 \geq y, \ldots, S_n \geq y)})$$

$$= \frac{1}{V(y)} \mathbb{E}_0(V(y+S_n) \mathbf{1}_{(A, S_1 \geq 0, \ldots, S_n \geq 0)}) = \frac{1}{V(y)} \mathbb{E}_0^+\left(\mathbf{1}_A \frac{V(S_n + y)}{V(S_n)}\right). \tag{A.3}$$

Notice that

$$\frac{V(S_n + y)}{V(S_n)} \leq 1 + \frac{V(y)}{V(S_n)} \leq 1 + V(y),$$

because the function $V(\cdot)$ is subadditive, increasing and $V(0) = 1$. Moreover, we have $V(S_n + y)/V(S_n) \to 1$ because $S_n \to \infty$, $\mathbb{P}^+$-a.s., hence we can apply dominated convergence when taking the $n \to \infty$ limit in (A.3) and (A.2) follows.

Observe that in particular we have proved that under $\mathbb{P}_z^+(\cdot | S_i \geq z\, \forall i \in \mathbb{N})$ the process $S - z$ has law $\mathbb{P}_0^+$.

For brevity we introduce the shorthand $S_{[a,b]}$ for the vector $(S_a, S_{a+1}, \ldots, S_b)$, and we write $S_{[a,b]} > x$ to mean $S_i > x$ for every $i = a, \ldots, b$. Then the pre-minimum and post-minimum processes may be expressed as $S_{[0,m]}$ and $S_{[m,\infty]} - S_m$ respectively. For $A, B \in \sigma(S_n, n \in \mathbb{N})$ we can write

$$\mathbb{P}_y^+(S_{[0,m]} \in A, S_{[m,\infty]} - S_m \in B) = \sum_{k \in \mathbb{N}} \int_{z \in [0,y]} \mathbb{P}_y^+(S_{[0,k]} \in A, S_{[m,\infty]} - S_m \in B, m = k, S_k \in \mathrm{d}z)$$

$$= \sum_{k \in \mathbb{N}} \int_{z \in [0,y]} \mathbb{P}_y^+(S_{[0,k]} \in A, S_{[0,k-1]} > z, S_k \in \mathrm{d}z) \mathbb{P}_z^+\left(B+z, \inf_{i \geq 0} S_i = z\right),$$

where we have used the Markov property of $\mathbb{P}_y^+$. Then applying (A.2) we obtain

$$\mathbb{P}_y^+(S_{[0,m]} \in A, S_{[m,\infty]} - S_m \in B) = \left(\sum_{k \in \mathbb{N}} \int_{z \in [0,y]} \mathbb{P}_y^+(S_{[0,k]} \in A, S_{[0,k-1]} > z, S_k \in \mathrm{d}z) \frac{1}{V(z)}\right) \mathbb{P}^+(B). \tag{A.4}$$

This factorization shows that under $\mathbb{P}_y^+$ the two processes $S_{[0,m]}$ and $S_{[m,\infty]} - S_m$ are indeed independent and the latter is distributed according to $\mathbb{P}^+$. It only remains to show that Eq. (A.1) holds true. For this observe that (A.4) yields in particular (just choose $B := \Omega_{\mathrm{RW}}$ and $A := \{S\colon S_m \geq x\}$)

$$\mathbb{P}_y^+(S_m \geq x) = \sum_{k\in\mathbb{N}} \int_{z\in[0,y]} \mathbb{P}_y^+(S_k \in [x,y], S_{[0,k-1]} > z, S_k \in \mathrm{d}z) \frac{1}{V(z)}$$

$$= \sum_{k\in\mathbb{N}} \int_{z\in[x,y]} \mathbb{P}_y^+(S_{[0,k-1]} > z, S_k \in \mathrm{d}z) \frac{1}{V(z)},$$

and by the definition (1.2) of $\mathbb{P}_y^+$ we get

$$\mathbb{P}_y^+(S_m \geq x) = \sum_{k\in\mathbb{N}} \int_{z\in[x,y]} \mathbb{P}_y(S_{[0,k-1]} > z, S_k \in \mathrm{d}z) \frac{V(z)}{V(y)} \frac{1}{V(z)}$$

$$= \frac{1}{V(y)} \sum_{k\in\mathbb{N}} \mathbb{P}_0(k \text{ is a ladder epoch}, S_k \in [x-y, 0]) = \frac{V(y-x)}{V(y)},$$

where we have used the definition (1.1) of the renewal function $V(\cdot)$.

We point out that one can give an explicit description of the pre-minimum process $\{S_k, k \leq m\}$. In fact this is closely related to the *random walk conditioned to die at zero* $(S, \mathbb{P}_z^{\searrow})$ described in Section 5, in analogy to the case of the Lévy process discussed in [11]. Let us work out the details in the lattice case, that is when the law of $S_1$ is supported in $\mathbb{Z}$ and is aperiodic.

Assume for simplicity that $y \in \mathbb{Z}^+$. We have already determined the law of the overall minimum $S_m$ under $\mathbb{P}_y^+$, namely

$$\mathbb{P}_y^+(S_m = x) = \frac{W(y-x)}{V(y)}, \quad x \in \{0, \ldots, y\}, \tag{A.5}$$

where the function $W(z) = V(z) - V(z-1)$ defined in (5.5) is the mass function of the renewal process $\{\overline{H}_k\}_k$ (we set $W(0) := V(0) = 1$ by definition).

Then to characterize the pre-minimum process it remains to give the joint law of the vector $(m, \{S_k - x\}_{0 \leq k \leq m})$ under $\mathbb{P}_y^+$ conditionally on $(S_m = x)$, for $x \in \{0, \ldots, y\}$. We claim that this is the same as the law of $(\zeta, \{S_k\}_{0 \leq k \leq \zeta})$ under $\mathbb{P}_{y-x}^{\searrow}$, where $\zeta$ denotes the first hitting time of zero.

Let us prove this claim. Notice that $x = y$ means that $m = 0$ and this squares with the fact that $\mathbb{P}_0^{\searrow}(\zeta = 0) = 1$. Therefore we may assume that $x \in \{0, \ldots, y-1\}$. We recall the notation $T_I := \inf\{n \in \mathbb{Z}^+ \colon S_n \in I\}$ for $I \subseteq \mathbb{R}$. Then for any $N \in \mathbb{N}$ and $A \in \sigma(S_1, \ldots, S_N)$, by (A.4) we can write:

$$\mathbb{P}_y^+(m > N, S_{[0,N]} - x \in A, S_m = x)$$

$$= \mathbb{P}_y(S_{[0,N]} - x \in A, T_{(-\infty, x]} > N, S_{T_{(-\infty, x]}} = x) \frac{1}{V(y)}$$

$$= \mathbb{P}_{y-x}(S_{[0,N]} \in A, T_{(-\infty, 0]} > N, S_{T_{(-\infty, 0]}} = 0) \frac{1}{V(y)}.$$

Next we apply the Markov property at time $N$, recalling that by definition $W(z) = \mathbb{P}_z(S_{T_{(-\infty, 0]}} = 0)$ for $z \in \mathbb{N}$, and using (A.5) we finally obtain

$$\mathbb{P}_y^+(m > N, S_{[0,N]} - x \in A | S_m = x)$$

$$= \frac{1}{W(y-x)} \mathbb{E}_{y-x}(S_{[0,N]} \in A, T_{(-\infty, 0]} > N, W(S_N))$$

$$= \mathbb{P}_{y-x}^{\searrow}(S_{[0,N]} \in A, \zeta > N),$$

where in the last equality we have applied (5.6).

*A.2. Proof of Lemma* 4.2

As in the proof of Lemma 4.1, it is convenient to rephrase the statement in terms of the unrescaled random walk. We recall that $m$ is the first instant at which $S$ reaches its overall minimum. We have to prove that if $(y_N)$ is a positive sequence such that $x_N := y_N/a_N \to 0$ as $N \to \infty$, then for every $\varepsilon > 0$

$$\lim_{N\to\infty} \mathbb{P}^+_{y_N}(m \geq \varepsilon N) = 0 \quad \text{and} \quad \lim_{N\to\infty} \mathbb{P}^+_{y_N}\Big(\sup_{0\leq k\leq m} S_k \geq \varepsilon a_N\Big) = 0.$$

We follow very closely the arguments in [7]. We have

$$\mathbb{P}^+_{y_N}(m > \varepsilon N) = \mathbb{P}^+_{y_N}\Big(\inf_{n\leq \lfloor\varepsilon N\rfloor} S_n > \inf_{n > \lfloor\varepsilon N\rfloor} S_n\Big)$$

$$= \int_{x\in[0,y_N]} \int_{z\in[x,\infty)} \mathbb{P}^+_{y_N}\Big(\inf_{n\leq \lfloor\varepsilon N\rfloor} S_n \in \mathrm{d}x, S_{\lfloor\varepsilon N\rfloor} \in \mathrm{d}z, \inf_{n>\lfloor\varepsilon N\rfloor} S_n < x\Big)$$

$$= \int_{x\in[0,y_N]} \int_{z\in[x,\infty)} \mathbb{P}^+_{y_N}\Big(\inf_{n\leq \lfloor\varepsilon N\rfloor} S_n \in \mathrm{d}x, S_{\lfloor\varepsilon N\rfloor} \in \mathrm{d}z\Big)\mathbb{P}^+_z\Big(\inf_{n\in\mathbb{N}} S_n < x\Big), \quad (A.6)$$

where in the last equality we have used the Markov property. Using the definition (1.2) of $\mathbb{P}^+$ and relation (A.1), we can rewrite the measure appearing in the integral above in terms of the unperturbed random walk measure $\mathbb{P}$: more precisely, for $x \in [0, y_N]$ and $z \geq x$ we obtain

$$\mathbb{P}^+_{y_N}\Big(\inf_{n\leq \lfloor\varepsilon N\rfloor} S_n \in \mathrm{d}x, S_{\lfloor\varepsilon N\rfloor} \in \mathrm{d}z\Big)\mathbb{P}^+_z\Big(\inf_{n\in\mathbb{N}} S_n < x\Big)$$

$$= \mathbb{P}_{y_N}\Big(\inf_{n\leq \lfloor\varepsilon N\rfloor} S_n \in \mathrm{d}x, S_{\lfloor\varepsilon N\rfloor} \in \mathrm{d}z\Big)\frac{V(z)}{V(y_N)}\frac{V(z)-V(z-x)}{V(z)}$$

$$\leq \mathbb{P}_{y_N}\Big(\inf_{n\leq \lfloor\varepsilon N\rfloor} S_n \in \mathrm{d}x, S_{\lfloor\varepsilon N\rfloor} \in \mathrm{d}z\Big),$$

where for the last inequality observe that $V(z) - V(z-x) \leq V(x) \leq V(y_N)$, because the renewal function $V(\cdot)$ is subadditive and increasing. Coming back to (A.6) we get

$$\mathbb{P}^+_{y_N}(m > \varepsilon N) \leq \int_{x\in[0,y_N]} \int_{z\in[x,\infty)} \mathbb{P}_{y_N}\Big(\inf_{n\leq \lfloor\varepsilon N\rfloor} S_n \in \mathrm{d}x, S_{\lfloor\varepsilon N\rfloor} \in \mathrm{d}z\Big)$$

$$= \mathbb{P}_{y_N}\Big(\inf_{n\leq \lfloor\varepsilon N\rfloor} S_n \in [0,y_N]\Big) \leq \mathbb{P}_0\Big(\inf_{n\leq \lfloor\varepsilon N\rfloor} S_n \geq -y_N\Big),$$

and since by hypothesis $y_N/a_N \to 0$, the last term above vanishes by the invariance principle for $\mathbb{P}_0$.

Next we pass to the analysis of the maximum. We introduce the stopping time $\tau_N := \inf\{k: S_k \geq \varepsilon a_N\}$. Taking $N$ sufficiently large so that $y_N/a_N \leq \varepsilon$, we have

$$\mathbb{P}^+_{y_N}\Big(\sup_{k\leq m} S_k \geq \varepsilon a_N\Big) = \mathbb{P}^+_{y_N}(\tau_N \leq m) = \mathbb{P}^+_{y_N}\Big(\inf_{k\leq \tau_N} S_k > \inf_{k>\tau_N} S_k\Big)$$

$$= \int_{x\in[0,y_N]} \int_{z\in[\varepsilon a_N,\infty)} \mathbb{P}^+_{y_N}\Big(\inf_{k\leq \tau_N} S_k \in \mathrm{d}x, S_{\tau_N} \in \mathrm{d}z\Big)\mathbb{P}^+_z\Big(\inf_{n\in\mathbb{N}} S_n < x\Big), \quad (A.7)$$

where we make use of the strong Markov property at $\tau_N$. Now it suffices to focus on the last factor: using relation (A.1) and the fact that $V(\cdot)$ is subadditive and increasing, for $x \leq y_N$ and $z \geq \varepsilon a_N$ we get

$$\mathbb{P}^+_z\Big(\inf_{n\in\mathbb{N}} S_n < x\Big) = \frac{V(z)-V(z-x)}{V(z)} \leq \frac{V(x)}{V(z)} \leq \frac{V(y_N)}{V(\varepsilon a_N)}.$$

Then plugging this into (A.7) we obtain simply

$$\mathbb{P}^+_{y_N}\left(\sup_{k\leq m} S_k \geq \varepsilon a_N\right) \leq \frac{V(y_N)}{V(\varepsilon a_N)} = \frac{V(x_N a_N)}{V(\varepsilon a_N)} \to 0, \qquad N \to \infty,$$

where the last convergence follows from the subadditivity of $V(\cdot)$ and from the fact that $x_N \to 0$.

## Appendix B. Conditioning to stay positive vs. nonnegative

We recall the definition of the event $\mathcal{C}_N := (S_1 \geq 0, \ldots, S_N \geq 0)$ and of the function $V(x) := \sum_{k\geq 0} \mathbf{P}(\overline{H}_k \leq x)$, where $\{\overline{H}_k\}_{k\geq 0}$ is the strict descending ladder heights process defined in the introduction. We also set $\mathcal{C}_N^\sim := (S_1 > 0, \ldots, S_N > 0)$ and we define a modified function $V^\sim(x) := V(x-) = \lim_{y\uparrow x} V(y)$ for $x > 0$, while for $x = 0$ we set $V^\sim(0) := \mathbf{E}(V^\sim(S_1)\mathbf{1}_{(S_1>0)})$. Then we have the following basic result.

**Proposition B.1.** *Assume that the random walk does not drift to* $-\infty$, *that is* $\limsup_k S_k = +\infty$, $\mathbb{P}$-*a.s. Then the function* $V(\cdot)$ *(resp.* $V^\sim(\cdot)$*) is invariant for the semigroup of the random walk killed when it first enters the negative half-line* $(-\infty, 0)$ *(resp. the nonpositive half-line* $(-\infty, 0]$*). More precisely one has*

$$V(x) = \mathbf{E}_x(V(S_N)\mathbf{1}_{\mathcal{C}_N}), \qquad V^\sim(x) = \mathbf{E}_x(V^\sim(S_N)\mathbf{1}_{\mathcal{C}_N^\sim}), \tag{B.1}$$

*for all* $x \geq 0$ *and* $N \in \mathbb{N}$.

**Proof.** Plainly, it is sufficient to show that (B.1) holds for $N = 1$, that is

$$V(x) = \mathbf{E}_x(V(S_1)\mathbf{1}_{(S_1\geq 0)}), \qquad V^\sim(x) = \mathbf{E}_x(V^\sim(S_1)\mathbf{1}_{(S_1>0)}), \tag{B.2}$$

and the general case will follow by the Markov property.

We first prove a particular case of (B.2), namely $V(0) = \mathbf{E}(V(S_1)\mathbf{1}_{(S_1\geq 0)})$, or equivalently

$$\int_{(y\geq 0)} \mathbf{P}(S_1 \in dy)V(y) = 1. \tag{B.3}$$

Setting $\widehat{S}_n := -S_n$, by the definition of $V(\cdot)$ we get

$$\int_{(y\geq 0)} \mathbf{P}(S_1 \in dy)V(y) = \int_{(y\geq 0)} \mathbf{P}(S_1 \in dy)\left(\sum_{k\geq 0}\sum_{n\geq 0} \mathbf{P}(\overline{T}_k = n, \widehat{S}_n \leq y)\right)$$

$$= \int_{(y\geq 0)} \mathbf{P}(S_1 \in dy)\left(\sum_{n\geq 0} \mathbf{P}(n \text{ is a ladder epoch}, \widehat{S}_n \leq y)\right)$$

$$= \mathbf{P}(S_1 \geq 0) + \int_{(y\geq 0)} \mathbf{P}(S_1 \in dy)\left(\sum_{n\geq 1} \mathbf{P}(\widehat{S}_1 > 0, \ldots, \widehat{S}_n > 0, \widehat{S}_n \leq y)\right),$$

where in the last equality we have applied the *Duality Lemma*, cf. [20], Chapter XII. Denoting by $T_1 := \inf\{n \geq 1 \colon S_n \geq 0\}$ the first *weak ascending* ladder epoch of $S$, we have

$$\int_{(y\geq 0)} \mathbf{P}(S_1 \in dy)V(y) = \mathbf{P}(T_1 = 1) + \int_{(y\geq 0)} \mathbf{P}(S_1 \in dy)\left(\sum_{n\geq 1} \mathbf{P}(S_1 < 0, \ldots, S_n \in [-y, 0))\right)$$

$$= \mathbf{P}(T_1 = 1) + \sum_{n\geq 1}\int_{(z<0)} \mathbf{P}(S_1 < 0, \ldots, S_n \in dz)\mathbf{P}(S_1 \geq -z)$$

$$= \mathbf{P}(T_1 = 1) + \sum_{n\geq 1} \mathbf{P}(S_1 < 0, \ldots, S_n < 0, S_{n+1} \geq 0) = \mathbf{P}(T_1 < \infty),$$

and since $\mathbf{P}(T_1 < \infty) = 1$, because by hypothesis $\limsup_k S_k = +\infty$, $\mathbb{P}$-a.s., Eq. (B.3) is proved.

Next we pass to the general case. Observe that $V(x) = \mathbf{E}(\mathcal{N}_{[0,x]})$ for $x \geq 0$ and $V^\sim(x) = \mathbf{E}(\mathcal{N}_{[0,x)})$ for $x > 0$, where for $I \subseteq \mathbb{R}^+$ we set $\mathcal{N}_I := \#\{k \geq 0 \colon \overline{H}_k \in I\}$. Then conditioning the variable $\mathcal{N}_{[0,x]}$ on $S_1$ and using the Markov property of $S$, we have for $x \geq 0$

$$V(x) = \int_{\mathbb{R}} \mathbf{P}(S_1 \in \mathrm{d}y)\{(1 + V(x+y) - V(y))\mathbf{1}_{(y \geq 0)} + (1 + V(x+y))\mathbf{1}_{(y \in [-x,0))} + \mathbf{1}_{(y < -x)}\}$$

$$= \mathbf{E}_x(V(S_1)\mathbf{1}_{(S_1 \geq 0)}) + 1 - \int_{(y \geq 0)} \mathbf{P}(S_1 \in \mathrm{d}y)V(y)$$

$$= \mathbf{E}_x(V(S_1)\mathbf{1}_{(S_1 \geq 0)}),$$

having used (B.3). Analogously we have for $x > 0$

$$V^\sim(x) = \int_{\mathbb{R}} \mathbf{P}(S_1 \in \mathrm{d}y)\{(1 + V^\sim(x+y) - V(y))\mathbf{1}_{(y \geq 0)} + (1 + V^\sim(x+y))\mathbf{1}_{(y \in (-x,0))} + \mathbf{1}_{(y \leq -x)}\}$$

$$= \mathbf{E}_x(V^\sim(S_1)\mathbf{1}_{(S_1 > 0)}) + 1 - \int_{(y \geq 0)} \mathbf{P}(S_1 \in \mathrm{d}y)V(y)$$

$$= \mathbf{E}_x(V^\sim(S_1)\mathbf{1}_{(S_1 > 0)}).$$

Finally, it suffices to observe that this relation holds true also for $x = 0$ by the very definition of $V^\sim(0)$, and the proof is completed. $\square$

## Acknowledgment

We thank Ron A. Doney for fruitful discussions.